\documentclass[11pt, letterpaper]{amsart}
\usepackage[utf8]{inputenc}

\usepackage{amssymb, amsmath, amsthm, wasysym, mathrsfs}
\usepackage{hyperref}
\usepackage{cleveref}
\usepackage[alphabetic,lite]{amsrefs}
\usepackage{tikz-cd}
\usepackage{comment}
\usepackage{thmtools}
\usepackage{setspace}
\usepackage{mathtools}
\usepackage{enumitem}
\usepackage{thm-restate}

\usepackage{caption}
\usepackage{subcaption}
\usepackage{float}

\newtheorem{lemma}{Lemma}[section]
\newtheorem{lemma*}{Lemma}
\newtheorem{theorem}[lemma]{Theorem}

\newtheorem{cor}[lemma]{Corollary}

\newtheorem{question}[lemma]{Question}
\newtheorem{claim*}{Claim}

\newtheorem{defn}[lemma]{Definition}

\theoremstyle{definition}
\newtheorem{remark}[lemma]{Remark}
\newtheorem*{lem}{Acknowledgements}

\theoremstyle{plain}

    \newtheoremstyle{TheoremNum}
        {\topsep}{\topsep}              
        {\itshape}                      
        {}                              
        {\bfseries}                     
        {.}                             
        { }                             
        {\thmname{#1}\thmnote{ \bfseries #3}}
    \theoremstyle{TheoremNum}
    \newtheorem{thmn}{Theorem}
    \newtheorem{corn}{Corollary}

\newcommand{\Hh}{{\mathbb H}}

\newcommand{\C}{{\mathbb C}}

\newcommand{\Z}{{\mathbb Z}}

\DeclareMathOperator{\N}{N}

\numberwithin{equation}{section}
\numberwithin{table}{section}
\title{Simply transitive geodesics and omnipotence of lattices in PSL$(2,\mathbb{C})$}
\author{Ian Agol, Tam Cheetham-West, Yair Minsky}
\date{Summer 2024}
\address{University of California, Berkeley, 970 Evans Hall \# 3840, Berkeley, CA 94720-3840}
\email{ianagol@math.berkeley.edu}
\address{Department of Mathematics \\ Yale University \\ New Haven, CT, 06511}
  \email{yair.minsky@yale.edu,tamunonye.cheetham-west@yale.edu}

\begin{document}
\pagestyle{plain}

\begin{abstract}
    We show that the isometry group of a finite-volume hyperbolic 3-manifold acts simply transitively on many of its closed geodesics. Combining this observation with the Virtual Special Theorems of the first author and Wise, we show that every non-arithmetic lattice in PSL$(2,\mathbb{C})$ is the full group of orientation-preserving isometries for some other lattice and that the orientation-preserving isometry group of a finite-volume hyperbolic 3-manifold acts non-trivially on the homology of some finite-sheeted cover.
\end{abstract}
\maketitle
\bibliographystyle{alpha}

\section{Introduction}

\medbreak\noindent
By the rigidity theorems of Mostow-Prasad, two finite-volume hyperbolic 3-manifolds $M_1,M_2$ with isomorphic fundamental groups are isometric. For finite-volume hyperbolic 3-manifolds, Bridson-Reid's conjecture (Conjecture 2.1 \cite{BridSurvey}) that finite covolume Kleinian groups are profinitely rigid (see Definition~\ref{prigid}) implies the profinite analog of Mostow-Prasad rigidity. Wilton-Zalesskii \cite{WZ1} showed that the profinite completion of a 3-manifold group detects whether the 3-manifold is hyperbolic of finite volume, Bridson-McReynolds-Reid-Spitler \cite{BMRS1} gave the first examples of profinitely rigid lattices in PSL$(2,\C)$, and Yi Liu \cite{Y} showed that the profinite completion distinguishes lattices in PSL$(2,\C)$ up to finite ambiguity. 
\medbreak\noindent The aforementioned theorems of Wilton-Zalesskii and Liu leverage specific properties of the fundamental groups of finite-volume hyperbolic 3-manifolds. In particular, it is crucial for the proofs of these theorems that the fundamental groups of finite-volume hyperbolic 3-manifolds are virtually special \cite{AgolHaken}\cite{WiseHaken}, that finite-volume hyperbolic 3-manifolds virtually fiber over $S^1$ \cite{Agolcriterion}, and that the fundamental groups of closed hyperbolic 3-manifolds have lots of virtually special quotients \cite{WiseHaken}. 
\medbreak\noindent In this note, we show that for non-arithmetic lattices in PSL$(2,\C)$ the general profinite rigidity question can be reduced to considering fundamental groups of fibered non-arithmetic hyperbolic 3-manifolds. To do this, we first show that there is a closed geodesic $\gamma\in M$ such that $f(\gamma)\ne\gamma$ for any non-trivial isometry $f:M\to M$. This is Lemma~\ref{freeorbitlemma}. By combining this lemma with the Malnormal Special Quotient Theorem of Wise \cite{WiseHaken}, and a general omnipotence theorem in the case of non-compact, finite-volume hyperbolic 3-manifold groups \cite{goodsam}, we show 
\medbreak
\begin{thmn}[\ref{mainthm} \& \ref{shepherd}]
  For $\Gamma_1$ a non-arithmetic torsion-free lattice in $\text{PSL}(2,\C)$, and $\Delta<\Gamma_1$ a finite-index subgroup, there is a lattice $\Gamma_2<\Delta<\Gamma_1$ such that $Aut^+(\Gamma_2)=\Gamma_1$.   
\end{thmn}
\medbreak\noindent
Here $Aut^+(G)$ denotes the orientation-preserving automorphisms of the lattice $G$. A similar theorem holds for the full automorphism group if we work over $Isom(\mathbb{H}^3)$.
\begin{defn}\label{prigid}
    A finitely generated residually finite group $\Gamma$ is {\it profinitely rigid} if for any finitely generated residually finite group $\Delta$ with the same set of finite quotients as $\Gamma$, $\Gamma\cong\Delta$. 
\end{defn} 
\medbreak\noindent As a corollary of these theorems, we have 
\medbreak\noindent
\begin{corn}[\ref{maincor}]
The following are equivalent:
    \begin{enumerate}
        \item The profinite rigidity of all non-arithmetic lattices in PSL$(2,\C)$. 
        \item The profinite rigidity of all fibered non-arithmetic lattices in PSL$(2,\C)$.
        \item The profinite rigidity of all special non-arithmetic lattices in PSL$(2,\C)$.      
    \end{enumerate}
\end{corn}
\begin{lem}
     {\it The authors thank Hee Oh, Alan Reid, Eduardo Reyes, and Franco Vargas Pallete for helpful comments and corrections. Ian Agol's research is supported by the Simons Foundation Award \#376200. Yair Minsky's research is supported by NSF Grant DMS-2005328.}
\end{lem}
    
\section{Omnipotence}
\begin{defn}
A finite subset $\{ g_1,\dots,g_n\}$ of infinite order elements in a group $G$ is {\bf independent} if the elements have pairwise non-conjugate non-trivial powers. 
\end{defn}
\begin{defn}
    A group $G$ is {\bf omnipotent} if for any independent subset $\{ g_1,\dots,g_n\}$ there is a constant $\kappa$ so that for any $n-tuple$ of positive integers $(e_1,\dots,e_n)$, there is a finite quotient $q:G\twoheadrightarrow Q$ with $ord(q(g_i))=\kappa e_i$ for $1\leq i\leq n$. 
\end{defn} 
\noindent The condition of omnipotence was first defined by Wise in \cite{omnipotence} where he proved that free groups are omnipotent (Theorem 3.5 \cite{omnipotence}). In \cite{WiseHaken}, Wise used the Malnormal Special Quotient Theorem (Theorem 12.2 \cite{WiseHaken}) to prove
\begin{theorem}[Theorem 14.26 \cite{WiseHaken}]
    Let $G$ be a virtually special, word-hyperbolic group. Then $G$ is omnipotent. 
\end{theorem}
\noindent Shepherd \cite{goodsam} proves a general omnipotence statement for finite independent subsets of {\it convex} (Definition 2.19 \cite{goodsam}) elements in virtually special cubulable groups. The following theorem is used in proving Theorems~\ref{shepherd} and \ref{moregeneral}.
\begin{theorem}[Theorem 1.2 \cite{goodsam}]\label{sammy}
    Let $G$ be a virtually special cubulable group. Then for any independent subset $\{ g_1,\dots,g_n\}$ of convex elements in $G$ there is a constant $\kappa$ so that for any $n-tuple$ of positive integers $(e_1,\dots,e_n)$, there is a finite quotient $q:G\twoheadrightarrow Q$ with $ord(q(g_i))=\kappa e_i$ for $1\leq i\leq n$. 
\end{theorem}

\section{Simply transitive geodesics}
\begin{theorem}[Theorem 1.1 (Corollary 2 \& 3)\cite{margthesis}, Proposition 5.4 \cite{gangw}]\label{marg}
Let $M$ be a finite-volume hyperbolic $n$-orbifold and let $$\pi^M_L=\#\{\text{oriented closed geodesics }\gamma\in M\,|\, l(\gamma)\leq L\}$$ Then $$\pi^M_L\sim \frac{e^{(n-1)L}}{(n-1)L}$$
\end{theorem}


\begin{lemma}\label{freeorbitlemma}
    Let $M$ be a finite-volume orientable hyperbolic $n$-manifold. There is a closed geodesic $\gamma\in M$ which is not fixed by any non-trivial isometry.
    \begin{proof}
  All geodesics in this proof will be unoriented. Let $X_l$ be the set of all closed geodesics of length $\leq l$ fixed by some non-trivial isometry of $M$. For $\gamma\in X_l$,
    \begin{enumerate}
        \item either $\gamma$ is fixed pointwise by $f$ and is a subset of Fix$(f)$, a proper, finite-area totally geodesic submanifold of $M$ (Type I) 
        \item $f$ is an involution on $\gamma$ fixing two points in $\gamma$ and interchanging the arcs between them and projecting to a geodesic arc of length $l/2$ perpendicular to the involution locus in $M/\langle f\rangle$ (Type II) or 
        \item $\gamma$ projects to a closed geodesic of fractional length $l/k$ in the orbifold $M/\langle f\rangle$ which we denote as $M_f$ (Type III). 
    \end{enumerate}
    Let $X_l^1$ be the set of all closed geodesics of length $\leq l$ in $\bigcup_{f\in Isom(M)}Fix(f)$. Let $X_l^2$ be the set of all Type II geodesics. Let $X_l^3$ be the set of all Type III geodesics. Now, $|X_l|=|X_l^1|+|X_l^2|+|X_l^3|\leq |X_l^1|+|X_l^2|+\sum_{f\in Isom(M)}\pi_{l/2}^{M_f}$. By Theorem~\ref{marg}, $$|X_l^1|\sim C_1\frac{e^{(n-2)l}}{(n-2)l},\enspace \pi_{l/2}^{M_f}\sim C_f\frac{e^{\frac{(n-1)l}{2}}}{l}$$ for some constants $C_1,C_f>0$. For any $f:M\to M$ an involution with fixed points, the involution locus of $M/\langle f\rangle$ is a proper non-empty properly immersed totally geodesic subset of $M/\langle f\rangle$, and Theorem 1.8 \cite{ohshah2}(see also Theorem 1 \cite{PPaulin}) gives an asymptotic for the number of common perpendiculars from the involution locus to itself. In particular, there is a constant $C_2>0$ such that $$|X_l^2|\sim C_{2,f}e^{\frac{(n-1)l}{2}}$$ Since $|X_l|\leq |X_l^1|+|X_l^2|+\sum_{f\in Isom(M)}\pi_{l/2}^{M_f}$, we have that $$|X_l|\leq C_1\frac{e^{(n-2)l}}{(n-2)l}+C_{2,f}e^{\frac{(n-1)l}{2}}+C_3\frac{e^{\frac{(n-1)l}{2}}}{l}\,(\clubsuit)$$ for a constant $C_3=\sum_{f\in Isom(M)}C_f$. As $l\to\infty$, the right-hand side of $(\clubsuit)$ $$ C_1\frac{e^{(n-2)l}}{(n-2)l}+C_{2,f}e^{\frac{(n-1)l}{2}}+C_3\frac{e^{\frac{(n-1)l}{2}}}{l}< \frac{e^{(n-1)l}}{(n-1)l}$$ Since $\pi^M_L\sim \frac{e^{(n-1)l}}{(n-1)l}$ by Theorem~\ref{marg}, it follows that $|X_l|< \pi_l^M/2$ as $l\to\infty$, and so there are many closed geodesics in $M$ that are not fixed by any isometry of $M$.
    \end{proof}
\end{lemma}
\begin{remark}
    In the literature, a special case of Lemma~\ref{freeorbitlemma} was proven by S. Kojima (Proposition 2 \cite{kojima}) for closed orientable hyperbolic 3-manifolds containing totally geodesic surfaces of genus $\geq 3$. 
\end{remark}
\begin{remark}
    If we restricted to simple closed geodesics, the statement of Lemma~\ref{freeorbitlemma} would be false. The hyperelliptic involution on the closed genus 2 surface fixes every simple closed geodesic (Theorem 1 \cite{haassusskind}). 
\end{remark}

\section{Main theorems for lattices in PSL$(2,\C)$}
\begin{theorem}\label{mainthm}
For $\Gamma_1$ a torsion-free cocompact non-arithmetic lattice in $\text{PSL}(2,\C)$, and $\Delta<\Gamma_1$ a finite-index subgroup, there is a lattice $\Gamma_2<\Delta<\Gamma_1$ such that $Aut^+(\Gamma_2)=\Gamma_1$. 
\begin{proof}
    For any finite index normal subgroup $\Delta<\Gamma_1$, $\Gamma_1<Aut(\Delta)$. Choose $\Delta'<\Delta<\Gamma_1$, a torsion-free finite-index subgroup fixed by the conjugation action of $\Lambda=Comm(\Gamma_1)$ the commensurator of $\Gamma_1$ in PSL$(2,\C)$ which is a lattice because $\Gamma_1$ is non-arithmetic (Theorem 1 p.2 \cite{margDisc}). 
    \medbreak\noindent By Lemma~\ref{freeorbitlemma}, there is a closed geodesic $\gamma$ in $\Hh^3/\Delta'$ on whose orbit the isometric action of $\Lambda/\Delta'$ is simply transitive. We further use $\gamma\in \Delta'$ to represent a choice of hyperbolic element whose conjugacy class corresponds to the geodesic $\gamma$. Set $g_1,\dots,g_k\in\Gamma_1/\Delta'<Out^+(\Delta')\cong \Lambda/\Delta'$ and $g_{k+1},\dots,g_n\in Out^+(\Delta')\setminus (\Gamma_1/\Delta')$, where $Out^+(G)$ denotes the orientation-preserving outer automorphisms. For each $g_i\in\Lambda/\Delta$ we choose $g_i'\in \Lambda$ a preimage of $g_i\in\Lambda/\Delta$. Since $\Delta'$ is hyperbolic and virtually special \cite{AgolHaken}\cite{WiseHaken}, we can apply the Malnormal Special Quotient Theorem (Theorem 12.2 \cite{WiseHaken}) to the independent collection of subgroups $\langle g_i'\gamma g_i'^{-1}\rangle<\Delta'$ to find a pair of non-zero integers $N_1\ne N_2$ such that $$\overline{\Delta'}\cong\Delta'/\langle (g_i'\gamma g_i'^{-1})^{N_1},1\leq i\leq k\,,(g_j'\gamma g_j'^{-1})^{N_2},k+1\leq j\leq n\rangle\,$$ is a hyperbolic, virtually special group. Moreover, the images of $(g_i'\gamma g_i'^{-1}), i\leq k$ and $(g_j'\gamma g_j'^{-1}), j>k$ will have orders $N_1$ and $N_2$ respectively by Theorem 7.2(1)  \cite{GrovesManning}. There is a natural action of $\Gamma_1$ on $\overline{\Delta'}$ and we can choose a finite-index subgroup $\overline{\Theta}<\overline{\Delta'}$ which is torsion-free and $\Gamma_1$-invariant. The preimage of $\overline{\Theta}$ in $\Delta'$ is a finite-index subgroup $\Theta$. Moreover, $$Aut^+(\Theta)=\Gamma_1$$ 
   To see this, it is sufficient to observe that $\Lambda$ is the full orientation-preserving automorphism group of $\Delta'$, and if any automorphism $f\in \Lambda\setminus\Gamma_1$ were to fix $\Theta$ as well, such $f$ would induce an automorphism $$\Tilde{f}:\Delta'/\Theta\to\Delta'/\Theta$$  
    which would send the image $\Tilde{g_i'}$ of a representative of a conjugacy class representing $\{g_i'\gamma g_i'^{-1}\}$ with $1\leq i\leq k$ to the image $\Tilde{g_j}'$ of a representative of a conjugacy class representing $\{g_j'\gamma g_j'^{-1}\}$ with $k+1\leq j\leq n$. The choice of $\overline{\Theta}$ torsion-free ensures that the torsion subgroups generated by $g_i'\gamma g_i'^{-1}$ for $1\leq i\leq n$ inject into $\Delta'/\Theta\cong \overline{\Delta'}/\overline{\Theta}$. By the construction of $\overline{\Delta'}$, the orders of $\Tilde{g_i}'$ and $\Tilde{g_j}'$ are distinct yielding a contradiction. 
\end{proof}
\end{theorem}

\begin{theorem}\label{shepherd}
For $\Gamma_1$ a torsion-free non-arithmetic non-uniform lattice in $\text{PSL}(2,\C)$ and any finite-index subgroup $\Delta<\Gamma_1$ there is a lattice $\Gamma_2<\Delta<\Gamma_1$ such that $Aut^+(\Gamma_2)=\Gamma_1$.  
\begin{proof}
As done in Theorem~\ref{mainthm}, for a finite index subgroup $\Delta<\Gamma_1$, we choose $\Delta'<\Delta<\Gamma_1$ a torsion-free, finite-index subgroup with $Aut(\Delta')=Comm(\Gamma_1)$ which is a lattice in PSL$(2,\C)$ because $\Gamma_1$ is non-arithmetic (Theorem 1 p.2 \cite{margDisc}). For brevity, we again denote $Comm(\Gamma_1)$ as $\Lambda$, and we choose a hyperbolic element $\gamma$ whose conjugacy class corresponds to a geodesic on which the full outer automorphism group of $\Delta'$ acts simply transitively. 
\medbreak\noindent For $g_1,\dots,g_k\in\Gamma_1/\Delta'$ and $g_{k+1},\dots,g_n\in (\Lambda/\Delta')\setminus (\Gamma_1/\Delta')$, we choose preimages $g_i'\in \Lambda$ for $1\leq i\leq n$. We consider the non-conjugate (in $\Delta'$) cyclic subgroups $<g_i'\gamma g_i'^{-1}>$ for $1\leq i\leq k$ and $<g_j'\gamma g_j'^{-1}>$ for $k+1\leq j\leq n$. By construction, these subgroups intersect all cusp subgroups of $\Delta'$ trivially. They are therefore convex subgroups of $\Delta'$ as defined in \cite{goodsam} (Definition 2.19). Since $\Delta'$ is cubulated and virtually special \cite{AgolHaken}\cite{WiseHaken}, by Theorem~\ref{sammy} $\Delta'$ is omnipotent and so there is an integer $\kappa$ and a finite quotient $\rho:\Delta'\rightarrow Q$ such that $ord(\rho(g_i'\gamma g_i'^{-1}))=\kappa N_1$ for $1\leq i\leq k$ and $ord(\rho(g_j'\gamma g_j'^{-1}))=\kappa N_2$ for $k+1\leq j\leq n$ and $N_1\ne N_2$. \medbreak\noindent Set $\Theta=\ker\rho$. For $f\in \Lambda\setminus\Gamma_1$, if $f(\Theta)=\Theta$, then $f$ induces a homomorphism $\Tilde{f}:Q\to Q$ that sends some element of order $\kappa N_1$ to an element of order $\kappa N_2$ which is impossible as $\kappa N_1\ne\kappa N_2$. Thus, $\Delta'<Aut(\Theta)<\Gamma_1$. If $\Theta$ is $\Gamma_1-$invariant, then $Aut(\Theta)=\Gamma_1$, and we are done. Otherwise, consider the $\Gamma_1-$invariant subgroup $\Theta'=\cap_{i=1}^k g_i\Theta g_i^{-1}$. 
\medbreak\noindent By construction, $\Delta'/\Theta'\cong \prod_{i=1}^k\Delta'/ g_i\Theta g_i^{-1}$ and we can check that the homomorphism $\rho':\Delta'\to \Delta'/\Theta'$ also satisfies $ord(\rho'(g_i\gamma g_i^{-1}))=\kappa N_1$ for $1\leq i\leq k$ and $ord(\rho'(g_j\gamma g_j^{-1}))=\kappa N_2$ for $k+1\leq j\leq n$ and $N_1\ne N_2$. To see that this is true, let $1\leq i,i'\leq k$ and let $k+1\leq j\leq n$. The order of the image of $g_{i'}\gamma g_{i'}^{-1}$ in $\Delta'/g_i\Theta g_i^{-1}$ is the same as the order of the image of $g_i^{-1}(g_{i'}\gamma g_{i'}^{-1})g_i$ in $\Delta'/\Theta$ which is $\kappa N_1$ by the previous paragraph. Thus, $ord(\rho'(g_{i'}\gamma g_{i'}^{-1}))$ is the least common multiple of the orders of the images of $g_{i'}\gamma g_{i'}^{-1}$ in each factor of $\prod_{i=1}^k\Delta'/ g_i\Theta g_i^{-1}$, and so $ord(\rho'(g_{i'}\gamma g_{i'}^{-1}))=\kappa N_1$. The order of the image of $g_{j}\gamma g_{j}^{-1}$ in $\Delta'/g_i\Theta g_i^{-1}$ is the same as the order of the image of $g_i^{-1}(g_{j}\gamma g_{j}^{-1})g_i$ in $\Delta'/\Theta$ which is $\kappa N_2$ by the previous paragraph. So, $ord(\rho'(g_{j}\gamma g_{j}^{-1}))$ is the least common multiple of the orders of the images of $g_{j}\gamma g_{j}^{-1}$ in each factor of $\prod_{i=1}^k\Delta'/ g_i\Theta g_i^{-1}$, and so $ord(\rho'(g_{j}\gamma g_{j}^{-1}))=\kappa N_2$ as claimed. Thus, $\Theta'$ is a $\Gamma_1$-invariant finite-index subgroup of $\Delta'$ with no other automorphisms outside $\Gamma_1$ and $Aut(\Theta')=\Gamma_1$ as claimed. 


\end{proof}
\end{theorem}

\medbreak\noindent We include a more general statement which applies to the Gromov-Piateski-Shapiro non-arithmetic hybrid lattices in SO$^+(n,1)$ \cite{GPS}. 
\begin{theorem}\label{moregeneral}
    Let $\Gamma<\,$SO$^+(n,1)$ be a virtually special non-arithmetic torsion-free lattice. Then for any finite index subgroup $\Delta<\Gamma$, there is a finite index subgroup $\Gamma_2<\Delta$ with $Aut(\Gamma_2)=\Gamma$.
    \begin{proof}
    The proof is the same as that of Theorem~\ref{shepherd}.
    \end{proof}
\end{theorem}
\begin{cor}\label{gps}
    Let $\Gamma<\,$SO$^+(n,1)$ be a cocompact non-arithmetic hybrid constructed in \cite{GPS}. Then for any finite index subgroup $\Delta<\Gamma$, there is a finite index subgroup $\Gamma_2<\Delta$ with $Aut(\Gamma_2)=\Gamma$.
    \begin{proof}
    By Proposition 9.1 \cite{BHW}, there is a finite-index subgroup $\Delta<\Gamma$ which is a quasiconvex subgroup of a simple type cocompact arithmetic lattice $\Gamma'<\,$SO$^+(n+1,1)$. The group $\Gamma'$ is virtually special by Theorem 1.6 \cite{HW1}, and therefore $\Gamma$ is virtually special by Proposition 7.2 \cite{HW2}. Thus, Theorem~\ref{moregeneral} applies to $\Gamma$. 
    \end{proof}
\end{cor}
\begin{remark}
    In the proofs of Theorems~\ref{mainthm} and \ref{shepherd}, it is crucial that the lattice is finite-index in its commensurator. In particular, this proof strategy does not apply to arithmetic lattices in PSL$(2,\C)$ and finitely generated groups with non-finitely generated (abstract) commensurators like (virtually) free groups. 
\end{remark}
\begin{remark}\label{BelL}
    A special case of Theorems~\ref{mainthm} and \ref{shepherd} for specially-defined non-arithmetic lattices in SO$^+(n,1)$ with epimomorphisms to (non-abelian) free groups was used by Belolipetsky-Lubotzky (see the proof of Theorem 3.1 \cite{BelLubotzkyIsom}) to show that for a fixed natural number $n\geq 2$, every finite group is the full isometry group of some finite-volume hyperbolic $n-$manifold (Theorem 1.1 \cite{BelLubotzkyIsom}). 
\end{remark}
\noindent Using Theorem~\ref{mainthm} and Theorem~\ref{shepherd}, we can prove:
\begin{cor}\label{maincor}
    The following are equivalent:
    \begin{enumerate}
        \item The profinite rigidity of all non-arithmetic lattices in PSL$(2,\C)$. 
        \item The profinite rigidity of all fibered non-arithmetic lattices in PSL$(2,\C)$.
        \item The profinite rigidity of all special non-arithmetic lattices in PSL$(2,\C)$.
        
    \end{enumerate}
    \begin{proof}
            It is sufficient to show that $(2)\implies (1)$ and $(3)\implies (2)$. To see that (2)$\implies$ (1), for a non-arithmetic lattice $\Gamma<\,$PSL$(2,\C)$ we apply Theorem~\ref{mainthm} and Theorem~\ref{shepherd} as follows; first, following the proofs of Theorem~\ref{mainthm} and Theorem~\ref{shepherd} and applying \cite{AgolHaken}\cite{WiseHaken}, we can choose $\Delta\triangleleft\Gamma$ such that $\Delta$ is the fundamental group of a hyperbolic 3-manifold that fibers over the circle. It will follow that the finite-index subgroup $\Theta<\Delta$ furnished by Theorem~\ref{mainthm} and Theorem~\ref{shepherd} will be a fibered lattice as well, satisfying $Aut^+(\Theta)=\Gamma$. By Theorem 4.4 \cite{BRPrasad}, once we assume that $\Theta$ is profinitely rigid, $\Gamma$ will be profinitely rigid as well. The proof that (3)$\implies$(2) follows along the same lines as the proof that (2)$\implies$(1). 
        \end{proof}
\end{cor}

\noindent We include an additional application of Lemma~\ref{freeorbitlemma}. 

\begin{theorem}\label{homology}
    Let $M$ be a finite-volume hyperbolic 3-manifold and let $f:M\to M$ be a non-trivial isometry. Then there is a finite-sheeted cover $M'\to M$ corresponding to a characteristic subgroup of $\pi_1(M)$ for which the induced isometry $f':M'\to M'$ acts non-trivially on $H_1(M';\Z)$.
    \begin{proof}
    By Lemma~\ref{freeorbitlemma}, there is a geodesic $\gamma$ with $f(\gamma)\ne\gamma$. After choosing an automorphism of $\pi_1(M)$ to represent $f$ and a hyperbolic element $\gamma\in\pi_1(M)$ to represent $\gamma$, we observe that for any characteristic finite-index subgroup $\pi_1(M')\triangleleft\pi_1(M)$, $f$ restricts to an isomorphism of $\pi_1(M')$ and for $n\in\N$ minimal such that $\gamma^n\in \pi_1(M')$ (also called the $\pi_1(M')$-{\it degree} of $\gamma$), $f(\gamma)^n\in\pi_1(M')$. Thus, $\gamma$ and $f(\gamma)$ have the same degree in every characteristic finite-index cover. Now, let $L=\langle \gamma\rangle$ and $T=\langle f(\gamma)\rangle$ be cyclic subgroups of $\pi_1(M)$. By construction $L,T$ are non-conjugate subgroups. 
    \medbreak\noindent When $M$ is compact, we apply the Malnormal Special Quotient Theorem (Theorem 12.3 \cite{WiseHaken}) to construct a word hyperbolic and virtually special quotient $\phi:\pi_1(M)\twoheadrightarrow G_L$ into which $L$ injects and where $T$ has finite image. By Lemma 14.12 \cite{WiseHaken}, there is a finite-index characteristic subgroup $J<G_L$ such that all conjugates of $\phi(L)$ have non-trivial images in the free abelianization of $J$. The preimage of $J$ in $\pi_1(M)$ is a finite-index subgroup $J'=\phi^{-1}(J)\triangleleft \pi_1(M)$ such that all conjugates of $L$ have non-trivial image in the free abelianization of $J'$ which is $H_1(J',\Z)_{\text{free}}$. 
    \medbreak\noindent Next, we pass to a characteristic finite-index subgroup $J''\triangleleft J'\triangleleft \pi_1(M)$, and let $n_{J''}$ be the $J''$-degree of $\gamma$. If the (non-trivial) homology class of $\gamma^{n_{J''}}$ in $H_1(J'',\Z)$ ($\cong H_1(\Hh^3/J'',\Z)$ since $\Hh^3/J''$ is aspherical) is the same as the homology class in $H_1(J'',\Z)$  of $f(\gamma^{n_{J''}})$, then the (non-trivial) homology class in $H_1(J',\Z)$ of $\gamma^{n_{J''}}$ is the same as the homology class in $H_1(J',\Z)$ of $f(\gamma^{n_{J''}})$. However, $f(\gamma^{n_{J''}})\in T$ which is distinguished from $L$ in the abelianization of $J'$ by the construction in the previous paragraph. Thus, in $H_1(J'',\Z)$, $f(\gamma^{n_{J''}})\ne \gamma^{n_{J''}}$, and that shows that the induced action on the homology of the characteristic finite-index cover corresponding to $J''$ is non-trivial.
    \medbreak\noindent In the case where $M$ is non-compact, we first choose a hyperbolic virtually special quotient $\rho:\pi_1(M)\to G$ for which $\rho(T)$ and $\rho(L)$ are non-conjugate infinite cyclic subgroups. We then apply the Malnormal Special Quotient Theorem to $G$ and the subgroups $\rho(T)$ and $\rho(L)$, just as in the compact case to find $G_L$ a hyperbolic, virtually special quotient of $G$ (and therefore $\pi_1(M)$) where $\rho(T)$ has finite image and $\rho(L)$ survives. For our choice of $G$, for example, we can set $G=\pi_1(\hat{M})$ where $\hat{M}$ is a compact hyperbolic 3-manifold  obtained from $M$ by hyperbolic Dehn filling, with $\rho:\pi_1(M)\to\pi_1(\hat{M})$ the Dehn filling epimorphism, such that $\rho(L)$ and $\rho(T)$ are non-conjugate infinite cyclic subgroups. One way to do this is to use strong conjugacy separability results (e.g. Theorem 1.1 \cite{Chagas2016Hyperbolic3G}) to choose a sufficiently large finite quotient $\pi_1(M)\twoheadrightarrow Q$, $|Q|<\infty$ where the images of $L,T,$ and the images of all cusp subgroups of $\pi_1(M)$ are non-trivial and pairwise non-conjugate. Since the images of all peripheral subgroups of $M$ are finite in $Q$, the finite quotient $\pi_1(M)\twoheadrightarrow Q$ factors through infinitely many Dehn fillings, and by Thurston's Dehn Surgery Theorem (\cite{Thurston} Theorem 5.8.2), we obtain $\hat{M}$ as required. By the Malnormal Special Quotient Theorem then, there is a quotient $G_L$ of $G$ with the specified properties (i.e. $G_L$ hyperbolic and virtually special, $L$ injects into $G_L$ and $T$ has finite image), the rest of the proof continues and concludes just as in the compact case above. 
    \end{proof}
\end{theorem}
\begin{remark}
    For a closed orientable surface $S$ the action of the isometry group is always faithful on $H_1(S,\Z)$ by a theorem of Hurwitz (see Theorem 6.8 \cite{farb2011primer}). In contrast, the isometry group of a hyperbolic 3-manifold can act homologically trivially (see the introduction of \cite{preni}). For example, when $M$ is a hyperbolic $\Z-$homology $3-$sphere (such as $1/n-$Dehn surgery on a hyperbolic knot in $S^3$ for sufficiently large $n$), $H_1(M;\Z)$ is trivial, and therefore, so is the action of $Isom^+(M)$ on $H_1(M,\Z)$. Theorem~\ref{homology} implies that there will be a finite-sheeted characteristic cover of $M$ for which $Isom^+(M)$ acts homologically faithfully. 
\end{remark}
\section{Remarks}

\noindent We conclude with some observations and questions coming from this circle of ideas. First, based on Remark~\ref{BelL} above,

\begin{question}
Is there a proof of Theorems~\ref{mainthm} and \ref{shepherd} using the methods of \cite{BelLubotzkyIsom} (lattice counting arguments)?
\end{question}

\begin{question}
    Is there an example of a non-maximal arithmetic lattice $\Gamma$ in PSL$(2,\C)$ for which Theorem~\ref{mainthm} holds?
\end{question}
\noindent For a maximal lattice $\Gamma$ in PSL$(2,\C)$, every normal subgroup $\Delta<\Gamma$ will have $Aut(\Delta)=\Gamma$. On the other hand, for any non-maximal arithmetic lattice $\Gamma<$ PSL$(2,\C)$, it would be remarkable if Theorem~\ref{mainthm} is true since by a theorem of Margulis, the commensurator of $\Gamma$ is dense in PSL$(2,\C)$ producing lots of hidden symmetries of $\Hh^3/\Gamma$. For any subgroup $\Delta<\Gamma$ with a hidden symmetry and a finite-index subgroup $\Theta$ characteristic in $\Delta$, $Aut(\Theta)\ne \Gamma$.

\begin{question}
Does Theorem~\ref{mainthm} hold for any non-maximal complex hyperbolic lattice? 
\end{question}
\noindent By the work of Stover, we know that there are pairs of non-isomorphic complex hyperbolic lattices with the same profinite completions. These pairs can be chosen to be commensurable (Corollary 1.3 \cite{stover2}) or non-commensurable (Theorem 1.1 \cite{stover1}). In the absence of a finiteness theorem such as Theorem 1 \cite{Y}, we can also ask about the complex hyperbolic analog of Theorem 4.4 \cite{BRPrasad}
\begin{question}
    If a lattice $\Gamma< \,$PU$(n,1)$ is profinitely rigid (among lattices in PU$(n,1)$), is its normalizer in PU$(n,1)$ profinitely rigid (among lattices in PU$(n,1)$)?
\end{question}
\noindent Finally, for a non-elementary hyperbolic group, we can ask whether an analog of Lemma~\ref{freeorbitlemma} holds.
\begin{question}
    Let $\Gamma$ be a non-elementary (relatively) hyperbolic group with $Out(\Gamma)$ non-trivial. For any non-trivial element $f\in Out(\Gamma)$, is there a (non-parabolic) primitive conjugacy class $\gamma\subset \Gamma$ for which $f(\gamma)\ne\gamma$?
\end{question}
\noindent One possible strategy towards a positive answer to this question is to use geodesic currents on groups i.e. to show that for a non-trivial outer automorphism $f\in Out(\Gamma)$, there is a geodesic current $\nu$ with $f_*(\nu)\ne \lambda\nu$ for $\lambda>0$ (where $f_*$ is the induced map on the space of currents) and then because conjugacy classes of elements approximate geodesic currents (Theorem 7 \cite{currents}), one may hope to argue that there will be a conjugacy class that is not fixed by the outer automorphism.


\bibliography{main}

\end{document}